\documentclass{article}
\usepackage{graphicx} 
\usepackage{amsthm}
\usepackage{fancyhdr} 
\usepackage{graphicx,tikz,xcolor} 

\title{Production optimization for agents of differing work rates}

\author{Peter M. Higgins}

\date{{}}

\begin{document}

\maketitle

\begin{abstract}
We devise a scheme for producing, in the least possible time, $n$
identical objects with $p$ agents that work at differing speeds.
This involves halting the process in order to transfer production
across agent types. For the case of two types of agent, we construct
a scheme based on the Euclidean algorithm that seeks to minimise the
number of pauses in production. 
\end{abstract}

\section{Introduction and context of the problem}

In {[}1{]} the author introduced the \emph{Biker-hiker problem}, which
entails finding optimal schemes for transporting $n$ travellers,
who collectively have $k$ bicycles, to their common destination in
the minimum possible time. Optimal schemes were those in which each
traveller rode $\frac{k}{n}$ of their journey by bicycle. Schemes
were represented by $n\times n$ binary matrices and an algorithm
was devised to determine when a matrix represented an optimal scheme,
based on the Dyck language of well-formed strings of parentheses.
A pair of mutually transpose matrices provided two optimal scheme
types, the first of which minimised the number of cycle handovers
while the transpose scheme kept the number and separation of the travelling
cohorts to a minimum. 

Here we re-interpret this as a problem of $n$ objects manufactured
by $n$ agents who have access to $k$ machines that work faster than
the agents. A modified version of the problem in which the machines
are more versatile in nature is the basis of this paper. We re-imagine
the scenario by saying that each of the $n$ agents has an identical
task to complete. For convenience of description, we take this task
to be the manufacture of an object. There are $k\leq n$ identical
machines available that can execute the task faster than an unsupported
agent. If a machine is not in use, an agent may continue to manufacture
the object they are constructing by passing it to that machine,
but only if the machine is configured to continue the build from this
point\emph{. }This constraint corresponds to the fact that in the
Biker-hiker problem a traveller may only mount a bicycle if they and
a bicycle are at the same point in the journey. In contrast, if an
agent takes over the manufacture of an object, they are capable of
recognising what point in the build has been reached for this object
and continue its construction from that point. This corresponds to
the ability of a traveller in the Biker-hiker context to continue
at any point of the journey on foot. 

To formulate the subject of this paper, we alter the nature of the
machines by allowing them to share with agents the capacity to continue
the build of an object from any given point in its construction. In
the original travelling setting, this would correspond to bicycles
that could instantaneously move from their current position to any
other in order to be used by a traveller.

Since the machines are faster than the agents, in any optimal scheme (one that minimises the total production time of the order), no machine will ever be idle. In view of this, there is no loss in adjusting the setting of the problem to view the $k$ machines as agents in their own right, and at all times some set of $k$ objects will be undergoing construction by the machines. There is then no need to have any more than $n-k$ other agents available to contribute to the building of the $n$ objects. 

Taken all together, these observations allow us to settle on the final
make-up of the problem. There are two types of agent, identical in
all respects except that one type works faster than the other. Let
us say they number $k_{1}$ and $k_{2}$, with $k_{1}+k_{2}=n$,
the total number of objects to be made. 

Indeed we shall widen the setting by allowing for an arbitrary number $m$ of agent types and for any number $p \geq n$ objects to be manufactured.  This general problem is formally stated in Section 2 where we present one optimal scheme type under the assumption that any time lost in halting the process in order to pass partially made objects between agents is negligible compared to the overall manufacturing time. In Section 3 however we return to the $m=2$ case and introduce methods that,
in general, greatly reduce the number of stoppages involved in executing an optimal scheme. 

\section{The general problem}

\subsection{Definition and principal features of an optimal scheme}

We have an order to manufacture $p\geq n$ identical objects using $n$ agents and we wish to do this in the minimum possible time. We shall refer to this challenge
as a \emph{$p$-object problem.} Each of the $n$ agents at our disposal produce one object at a time, but with varying speeds.
There are $k_{i}$ agents of type $i$ and $m$ types are available,
so that $k_{1}+k_{2}+\cdots+k_{m}=n$. Agent type $i$ takes $t_{i}$
time units (which we shall call hours) to complete the manufacture
of one object. A stipulated process for completing the order will be known as an \emph{$p$-scheme}, or simply a \emph{scheme}. 

At any point during the process an agent may be halted and its partially
constructed object replaced by another partially built object. This
second object may be at a different stage of construction but all
agents are capable of recognising this and can continue the build
from the current state. Our initial analysis will assume that the
time required to pass a partially built object from one agent to another
is negligible compared to the overall build time. 

Observe that if a scheme $S$ has the property that no agent is ever idle and all of them simultaneously complete the build of the object in hand, then $S$ is  optimal for such a scheme is working at maximum capacity all throughout the execution of the order. Once we exhibit the existence  of such a scheme $S$ in the general case, which we do in Section 2.2, it will follow  that optimal schemes are characterised by this property of simultaneous completion, for any scheme that lacks it will not be working to full capacity at some point during the build. We illustrate this principle through a simple example.

\vspace{.2cm}

\textbf{Example 2.1.1 } Take $k_{1}=k_{2}=1$ (so that $n=m=2$), and let
$t_{1}=1$, $t_{2}=2$. In order for both agents to work continuously
until the order is completed, we exchange objects at just the right
moment, which in this case is the $40$ minute mark, and then both
objects are fully built after $80$ minutes. 

\vspace{.2cm}

\textbf{Proposition 2.1.2 } Given the existence of an optimal solution scheme $S$  for the $n$-object problem (featuring $n$ agents), an optimal solution exists for the $p$-object problem for any $p\geq n$.

\begin{proof} 

 Write $p = an+b$ for a positive integer $a$ and a non-negative integer $b$, with $b$ satisfying $0\leq b\leq n-1$. Act the given optimal scheme $S$ for the $n$-object problem to produce the first set of $n$ manufactured objects.  Then repeat $S$ a further $a - 1$ times, yielding 
an output of $an$ objects, which have been manufactured in the least possible time as no agent has had an idle moment. There remain $b$ objects still to be produced, and so we act the given optimal $b$-object scheme on a set of $b$ of our fastest agents to complete the process in minimum time, thereby solving our optimization problem. 

\end{proof}

It follows that in order to find examples of optimal schemes we may henceforth restrict attention to the case where $p = n$, and so the numbers of agents and objects match. 

\vspace{.2cm}
 
{\bf Proposition 2.1.3} (a) In any optimal scheme $S$ for the $n$-object problem:

(i) the time $H$ required to complete the scheme is the harmonic
mean of the completion times of the individual agents:
\begin{equation}
H=n(\sum_{i=1}^{m}k_{i}t_{i}^{-1})^{-1};
\end{equation}

(ii) the proportion $p_{j}$ of the build constructed by the set of
type $j$ agents is: 
\begin{equation}
p_{j}=\frac{k_{j}}{t_{j}}(\sum_{i=1}^{m}k_{i}t_{i}^{-1})^{-1}.
\end{equation}

(b) Conversely, if all agents of a scheme begin simultaneously and work for duration $H$ as specified by (1),  then in doing so they have completed an optimal scheme.

\begin{proof}
    
(a)(i) In an optimal solution all agents work at full
capacity for some common time length, $H$. Since an agent of type
$i$ takes $t_{i}$ hours to make an object, its production rate is
$t_{i}^{-1}$ objects/hour. The combined rate of production, $R$
of all the agents in objects/hour is therefore given by the sum:
\[
R=\sum_{i=1}^{m}k_{i}t_{i}^{-1}.
\]

The time taken for the agents to produce the equivalent of $1$ object
is therefore $R^{-1}$, and so the total time $H$ to manufacture
the order of $n$ objects is given by $nR^{-1}$, which is the harmonic
mean of the individual times, as stated in (1).

\vspace{.2cm}

(ii)\textbf{ }The production rate of any agent of type $j$ is $t_{j}^{-1}$
objects/hour. Collectively, while executing an optimal scheme, the
$k_{j}$ agents of type $j$ produce $\frac{k_{j}H}{t_{j}}=\frac{k_{j}n}{t_{j}}(\sum_{i=1}^{m}k_{i}t_{i}^{-1})^{-1}$
objects, and so the proportion $p_j$ of the $n$ objects produced by the
type $j$ agents is as stated in (2):
\[
p_{j}=\frac{k_{j}}{t_{j}}(\sum_{i=1}^{m}k_{i}t_{i}^{-1})^{-1}=\frac{k_{j}}{t_{j}R}.
\]

(b) The combined work rate of the set of agents is $R$, and since $H=nR^{-1}$ it follows that after time $H$ the agents have produced $HR = n$ objects.  Therefore the production is complete and their action represents an optimal scheme. 

\end{proof}

The partition $P$ of the time interval $I$ of duration $H$  into $n$ equal intervals is the basis of the fundamental optimal scheme we shall introduce, as it is the length of time for agents of an optimal scheme to collectively build the equivalent of one object.

\vspace{.2cm}

\textbf{Definition 2.1.4 }(a) We shall refer to $H$ as defined in (1) as the \emph{harmonic optimum} time for an optimal $n$-object scheme. 

(b) The partition $P$ divides $I$ into $n$ intervals each of length $\frac{H}{n}=R^{-1}$. We call $R^{-1}$ the 
\emph{atomic time unit }(a.u.) for an optimal scheme.

(c) A scheme is \emph{uniform} if, for all $i$ $(1\leq i\leq m)$, each object is worked continuously and is worked by type $i$ agents for exactly $k_{i}$ a.u. We shall refer to $k_{i}$ as the \emph{type $i$ quota }for
objects in a uniform scheme. 

\vspace{.2cm}

\textbf{Proposition 2.1.5 }(a) Let $S$ be an optimal scheme for the
$n$-object problem. If, for all $i$, the time each object is worked
by type $i$ agents is the same for all objects, then $S$ is a uniform
scheme. 

(b) All uniform schemes are optimal.

(c) For $m=2$, a scheme $S$ is optimal if and only if $S$ is uniform. 

\begin{proof}
    
(a) Since $S$ is optimal, all agents work for $n$
a.u., and so the sum of the total times worked by type $i$ agents
is $nk_{i}$ a.u. If each of the $n$ objects were worked by type
$i$ agents for $K_{i}$ a.u., then the sum total work time by the
type $i$ agents would also equal $nK_{i}$, whence $K_{i}=k_{i}$,
and so $S$ is uniform. 

(b) Suppose that $S$ is a uniform scheme for the $n$-object problem.
Then each object is worked for $k_{j}$ a.u. by type $j$ agents.
Since one a.u. is equal to $\big(\sum_{i=1}^{m}k_{i}t_{i}^{-1}\big)^{-1}$,
the proportion of each object produced in $H=n$ a.u. is:
\[
\sum_{j=1}^{m}k_{j}t_{j}^{-1}\big(\sum_{i=1}^{m}k_{i}t_{i}^{-1}\big)^{-1}=\sum_{j=1}^{m}p_{j}=1.
\]
Therefore each object is completed after time $H$, and so $S$ is
optimal. 

(c) If $S$ is uniform, then $S$ is optimal by (b) so it only remains
to check that if $m=2$ and $S$ is an optimal scheme then $S$ is
uniform. For any pair of objects $O$ and $O'$, let $K_{1}$ and
$K_{1}'$ be the respective times, measured in a.u. that each object
is worked by type $1$ agents. Then since any two objects $O$ and
$O'$ have equal times of manufacture in an optimal scheme, we have:
\[
K_1t_1^{-1}+(n-K_1)t_2^{-1}=K_{1}'t_1^{-1}+(n-K_{1}')t_2^{-1}
\]
\[
\Leftrightarrow K_1(t_1^{-1} - t_2^{-1})=K_{1}'(t_1^{-1}-t_2^{-1}). \]
Since $t_1\neq t_2$, we get $K_1=K_{1}'$, whence it follows from (a) that $S$ is uniform.

\end{proof}

\textbf{Example 2.1.6 } For $m\geq3$ it is not necessarily the case
that all optimal schemes are uniform. To see this take $k_{1}=k_{2}=k_{3}=1,$
so that $n=m=3$. Let the completion times of agents $A_{1},A_{2},$
and $A_{3}$ be respectively $3,6,$ and $4$. 
\[
H=3(\sum_{i=1}^{3}k_{i}t_{i}^{-1})^{-1}=3(\frac{1}{3}+\frac{1}{6}+\frac{1}{4})^{-1}=3(\frac{4+2+3}{12})^{-1}=4.
\]
There is a non-uniform optimal scheme $S$ constructed as follows.
Agents $1$ and $2$ exchange objects $O_{1}$ and $O_{2}$ after
$2$ hours while agent $3$ works $O_{3}$ for $4$ hours until completion.
After $4$ hours the proportion of objects $1$ and $2$ that have
been completed will be, in both cases, $2(\frac{1}{3}+\frac{1}{6})=1$,
and so all three objects are completed in $H=4$ hours. Therefore
$S$ is optimal but $S$ is not uniform as it is not the case that
each object is worked by each type of agent for the same length of
time. We shall expand on this in Section 2.2.

\subsection{Optimal scheme construction}

\textbf{2.2.1 }\emph{The $n$-cyclic scheme.}

\vspace{.2cm} 

We label the $n$ agents $A_{1},A_{2},\cdots,A_{n}$, a freely chosen
order. We label the $n$ objects $O_{1},O_{2},\cdots,O_{n}$, where
$O_{i}$ is the object that begins to be made by agent $A_{i}$ in
the first interval $I_{1}$. The objects are then cycled around through
the agents: picture all the agents arranged in a circle. At the end
of each time interval $I_{j}$, a whistle blows. Each agent then takes
the object it is currently working on and passes it to the agent on
its right to continue the build. In symbols, at the end of interval
$I_{j}$, the object currently being worked on by $A_{i}$ is handed
to $A_{i+1}$, where we take $n+1=1$. 

\newpage

\textbf{Theorem 2.2.2 }The $n$-cyclic scheme, $C_{n}$, is optimal. 

\begin{proof}
    
By examining the evolution of $O_{i}$ in $C_{n}$
we see that in the successive intervals $I_{1},I_{2},\cdots,I_{j},\cdots,I_{n}$
object $O_{i}$ is worked by the respective agents, 
\[
A_{i},A_{i+1},\cdots,A_{i+j-1},\cdots,A_{i+n-1}=A_{i-1}.
\]
Therefore $O_{i}$ is worked on by each of the $n$ agents exactly once, and for the same length of time, which is $1$ a.u. Since $H = n$ a.u., it follows from Proposition 2.1.3(b) that $S$ is optimal.  

\end{proof}

\textbf{Example 2.2.4 }Let us take $k_{1}=3,k_{2}=4$, and $k_{3}=1$,
(so that $n=8$), with $t_{1}=1,t_{2}=2$, $t_{3}=4$. Then in objects/hour we have
\[
R=\sum_{i=1}^{3}k_{i}t_{i}^{-1}=\frac{3}{1}+\frac{4}{2}+\frac{1}{4}=\frac{21}{4},  H=\frac{n}{R}=\frac{32}{21}=1\frac{11}{21}.
\]
Each of the $8$ intervals equals $\frac{H}{n}=\frac{1}{R}=\frac{4}{21}$ hours $=\frac{4}{21}\times60=\frac{80}{7}=11\frac{3}{7}$
minutes. The proportion of the objects built by each agent type is:
\[
p_{1}=\frac{k_{1}}{Rt_{1}}=\frac{3\times4}{21\times1}=\frac{12}{21}=\frac{4}{7},\,p_{2}=\frac{4\times4}{21\times2}=\frac{8}{21},\,p_{3}=\frac{1\times4}{21\times4}=\frac{1}{21}.
\]

Let us choose to order the agents in increasing execution time. We
may track the build of any of the objects. For instance at the end
of the sixth stage, $O_{5}$ will have been worked on by type 2 agents
for the first three stages, ($A_{5},A_{6},$ and $A_{7}$, representing
$\frac{3}{4}$ of the type $2$ quota), for the fourth stage by the
type 3 agent $A_{8}$, and for the fifth and sixth stages by the type
1 agents, $A_{1}$ and $A_{2}$ (representing $\frac{2}{3}$ of the
type 1 quota). Hence the proportion of the build of $O_{5}$ completed
at the end of Stage 6 is: 
\[
\frac{3}{4}p_{2}+p_{3}+\frac{2}{3}p_{1}=\frac{3}{4}\times\frac{8}{21}+\frac{1}{21}+\frac{2}{3}\times\frac{12}{21}=\frac{15}{21}=\frac{5}{7}\approx71.4\%.
\]

\subsection{Schemes that run uniform sub-schemes in parallel}

Example 2.1.6 was a simple instance of two uniform schemes running in parallel. We explore this idea further using properties of harmonic means. 

\vspace{.2cm}

\textbf{Lemma 2.3.1 }(a) Let $T$ be a finite list of positive numbers
(so that repeats are allowed) with harmonic mean $H(T)=H$. 
Let $U$ be a sub-list of $T$ strictly contained in $T$. 
If $H(U)=H$ then $H(T\setminus U)=H(U)=H$. 

(b) If $U,V$ are finite lists of positive numbers with $H(U)=H(V)=h$, then $H(U\cup V)=h$, where $U\cup V$ denotes the union of $U$
and $V$ as lists, so there is a distinct element of $U\cup V$ for
each member of $U$ and of $V$. 

\begin{proof}
   (a) By re-indexing as necessary, we may write $U=\{a_{1},a_{2},\cdots,a_{m}\}$
and $T\setminus U=\{a_{m+1},a_{m+2},\cdots,a_{m+n}\}$ for some $m,n\geq1$.
Let $S_{1}=\sum_{i=1}^{m}a_{i}^{-1}$, $S_{2}=\sum_{i=m+1}^{m+n}a_{i}^{-1}$,
and $S=S_{1}+S_{2}$. We are given that
\[
\frac{S}{m+n}=\frac{S_{1}}{m}=\frac{1}{H}
\]
\[
\Rightarrow S_{2}=S-S_{1}=\frac{m+n}{H}-\frac{m}{H}=\frac{n}{H}.
\]
\[
\Rightarrow H(T\setminus U)=\frac{n}{S_{2}}=n\cdot\frac{H}{n}=H.
\]

(b) We may take $U=\{a_{1},a_{2},\cdots,a_{m}\}$ and $V=\{a_{m+1},a_{m+2},\cdots,a_{m+k}\}$,
for some $k\geq1$ and let $S_{1}$ and $S_{2}$ denote the respective
sums of the reciprocals of the members of $U$ and of $V$. We are
given
\[
\frac{S_{1}}{m}=\frac{S_{2}}{k}=\frac{1}{h}
\]
\[
\Rightarrow S_{1}+S_{2}=\frac{m+k}{h}
\]
\[
\Rightarrow  H(U\cup V)=\frac{m+k}{S_{1}+S_{2}}= (m + k)\cdot \frac {h}{m + k} = h. 
\]
\end{proof}

\textbf{Remark 2.3.2} The conclusion of Lemma 2.3.1(b) is not true if we take the set union of $U$ and $V$ when they have elements in common. For example, take $U=\{1,5,20\},\,V=\{1,6,12\}$ so that, as sets, $U\cup V=\{1,5,6,12,20\}$.
Then $H(U)=H(V)=\frac{12}{5}$ but $H(U\cup V)=\frac{10}{3}$. However, if we treat the collections as lists, allowing repeated members when 
taking the union, (in this example this would yield the list $(1,1,5,6,12,20)$),
then the argument of the lemma applies and all three harmonic means
coincide. 

\vspace{.2cm}

Let $T=\{t_{1},t_{2},\cdots,t_{n}\}$, a list of job completion times
of our $n$ agents in the $n$-object problem, and let $H$ be the
harmonic mean of $T$. If $T$ contains a proper sub-list $U$ with
$H(U)=H$, then by Lemma 2.3.1(a), we have $H(T\setminus U)=H$ also,
allowing us to split $T$ into two complementary sub-lists $U,V$
with $H(U)=H(V)=H$. This process may be repeated on the lists $U$
and $V$ until we have partitioned $T$ into a disjoint union of $r$ lists say, which we write as 

$$T=T_{1}\bigoplus T_{2}\bigoplus\cdots\bigoplus T_{r}\,\,(r\geq1),$$ where each sub-list $T_{i}$ is \emph{irreducible}, meaning that $T_{i}$
has no proper sub-list $U$ with $H(U)=H$. We call this an \emph{irreducible} representation of $T$. 

Given an irreducible representation of $T$, we may construct an optimal scheme $S$ for the $n$-object problem from any collection of optimal schemes $\{S_{i}\}_{1\leq i\leq r}$, where $S_{i}$ is an optimal scheme for the $n_{i}$-object problem, with $n_{i}$ defined by $n_{i}=\sum_{j=1}^{s_{i}}t_{i,j},$ where $t_{i,1},t_{i,2},\cdots,t_{i,s_{i}}$
are the members of the list $T$ that belong to $T_{i}$. Since the harmonic means of all these $r$ sub-problems equal $H$, we may run these $r$ schemes $S_{i}$ in parallel to provide an optimal solution $S$ of the original $n$-object problem. In recognition we denote this scheme as 
\[
S=S_{1}\bigoplus S_{2}\bigoplus\cdots\bigoplus S_{r}.
\]

In our Example 2.1.6, we have $T=\{3,6,4\}$, so that $T_{1}=\{3,6\},$
$T_{2}=\{4\}$, and $H=4$. The sub-schemes $S_{1}$ and $S_{2}$
are respectively the $2$-cyclic and $1$-cyclic schemes. The overall
optimal scheme $S$ is then the sum of two irreducible uniform schemes,
$S=S_{1}\bigoplus S_{2}$, but $S$ is not itself uniform. 

A given finite list $T$ of positive numbers may always be presented as a sum  of irreducible harmonic components. However, as our next example shows,  $T$ may have more than one irreducible representation. \footnote{This example was devised by Alexei Vernitski.}

\vspace{.2cm}

\textbf{Example 2.3.3 }The set $T=\{2,3,4,5,6,7,9,10,12,14,15\}$
can be split into two subset pairs preserving the harmonic mean
either as $U_1=\{2,7,9,10,15\}$ and $U_2 = T\setminus U_1 = \{3,4,5,6,12,14\}$,
or as $V_1 = \{2,5,6,10,14,15\}$ and $V_2=T\setminus V_1 = \{3,4,7,9,12\}$,
but not in any other ways. All five harmonic means come to $\frac{315}{58}\approx5.4310$. 

\vspace{.2cm}

We would expect that the problem of determining whether a given list admits a harmonic partition is hard as the  simpler problem of determining whether such a list has an additive partition is NP-complete, [7]. However we show below how to generate any number of such examples of lists with multiple harmonic partitions.  This allows us to construct an optimal non-uniform scheme that is not a parallel sum of irreducible schemes. 

\vspace{.2cm}

\textbf{Proposition 2.3.4} The harmonic mean $H=H(x,y)$ of two positive real numbers $x$ and $y$ satisfies $H(x,y)=2m$ if and only if $(x-m)(y-m)=m^{2}$. 

\begin{proof}
    
\[
\big(H(x,y)=2m\big)\Leftrightarrow\big(\frac{2xy}{x+y}=2m\big)\Leftrightarrow\big(xy-m(x+y)=0\big)
\]
\[
\Leftrightarrow\big(xy-mx-my+m^{2}=m^{2}\big)\Leftrightarrow\big((x-m)(y-m)=m^{2}\big).
\]

\end{proof}

This result lets us produce any number of distinct pairs of integers, $x$ and $y$, with a common harmonic mean through choosing $m$ so that $m^{2}$ has the requisite number of pairs of distinct factors. The
structure of our next example requires three such pairs, so we take $m  = 6$ and employ the factorizations $36=2\times18=3\times12=4\times9$. 

\vspace{.2cm}

\textbf{Example 2.3.5 } We construct a scheme $S$ with $n=6$ agents and objects and with respective job completion times of these agents
$A_1,\cdots,A_6$ given by $t_1 = 2 + 6,t_2 = 18+6$, so that $t_1 = 8,t_2 = 24$;
similarly $t_3 = 9,t_4 = 18$, and $t_5 = 10,t_6 = 15$. By Proposition
2.3.4, each of these pairs share a harmonic mean of $H=2\times6=12$.
Next, it follows by Lemma 2.3.1(b) that the harmonic mean of the union
of any two of the pairs of agent times, $(8,24),(9,18)$ and $(10,15)$, is also $12$, whence it follows again by Lemma 2.3.1(b) that $12$ is likewise the harmonic mean of the full set of six agent production times. Consider the scheme $S$ represented by the table in Figure 1. 

\vspace{.2cm}

~~~~~~~~~~~~~~~~~~~~~~~~~~~~~~%
\begin{tabular}{|c|c|c|c|c|c|}
\hline 
$O_{1}$ & $O_{2}$ & $O_{3}$ & $O_{4}$ & $O_{5}$ & $O_{6}$\tabularnewline
\hline 
\hline 
1 & 2 & 3 & 4 & 5 & 6\tabularnewline
\hline 
2 & 1 & 4 & 3 & 6 & 5\tabularnewline
\hline 
3 & 4 & 5 & 6 & 1 & 2\tabularnewline
\hline 
4 & 3 & 6 & 5 & 2 & 1\tabularnewline
\hline 
\end{tabular}

\vspace{.2cm}

\centerline{Figure 1: table for an optimal non-uniform scheme.}

\vspace{.2cm}

Each object $O_i$ is worked successively by the agents in its column when passing from top to bottom. For example $O_4$ is worked in turn by agents $A_4,A_3,A_6$, and $A_5$. Since each object
is worked by the agents of two pairs, each with the harmonic mean of $12$ hours, it follows that the duration of the common interval between exchanges is $\frac{12}{4}=3$ hours. Each agent appears exactly once in each of the $3$-hour windows, which are represented by the rows. It follows that $S$ is indeed an optimal scheme. However each pair of objects share common agents and so $S$ is not a sum of irreducible schemes, and nor is $S$ uniform as each object is worked by only $4$ of the $6$ agents. 

As a bonus we note that since it consists of three pairs with equal harmonic mean of $12$, the set  $\{8,9,10,15,18,24\}$, can be split into three harmonic partitions, one for each pair.  We can indeed create an example of order five by taking the set consisting of two of these pairs and their common harmonic mean of $12$: $T=\{8,9,12,18,24\}$, the two partitions being 
\[U = \{8,24\}, T\setminus U = \{9,12,18\}  \hspace{.2cm} \& \hspace{.2cm}
V =  \{9,18\},T\setminus V = \{8,12,24\}. \]

\subsection{Removing unnecessary exchanges}

It is always the case that we may take $m\leq n$. Indeed in most
real world examples we will have $m<<n$ for typically $n$ might
be of the order of hundreds or thousands, while $m$, the number of
different manufacturing speeds of our agents, might be in single digits.
And so even though we were operating only a handful of different agent
types, the cyclic scheme algorithm could involve thousands of stoppages
of the manufacturing process. We therefore look to see if we can adjust
our scheme to lower the number of stoppages. For a given scheme of
production $S$, we shall call the number of times the whistle blows to halt production the \emph{halting number}, and denote it by $h = h(S)$.
For the $n$-cyclic scheme $C_n$ we have $h(C_n)=n-1$. 

\vspace{.2cm}

\textbf{Remark 2.4.1 } Let $n'=\frac{n}{d}$,
where $d$ is the greatest common divisor of the list of integers
$k_1,k_2,\cdots,k_m$. Note that $d$ is also a divisor of $n$. Collect the $n$ objects into $n'$ groups of $d$ objects, regarding each such $d$-set as a single object. In a similar way group the $k_{i}$ agents of type $i$ into $k_{i}'=\frac{k_{i}}{d}$ sets of $d$ agents, with each set regarded as a single agent. We may now view this $n$-object problem as an 
$n'$-object problem with a list of $k_{1}',k_{2}',\cdots,k_{m}'$ agents. In particular, this procedure would allow replacement of a $C_n$-scheme by a $C_{n'}$-scheme, thereby reducing the halt number by a factor of $d$.   A typical object of the reduced scheme will consist of $d$ objects, all worked to the same point in their manufacture, which is then acted on in parallel by a collection of $d$ agents
of the same type to continue the scheme. 

Henceforth we will assume that the reduction process described has been carried out so that the gcd $d$ of $k_1,k_2,\cdots,k_m$ is $1$. 

\section{The $m = 2$ case: the Euclidean scheme}

The first interesting value for $m$ is therefore $m = 2$, which is
the focus of the remainder of the paper. By Proposition 2.1.5(c),
the only optimal schemes are uniform. For $n=2$, $C_{2}$ is the
unique optimal scheme with the fewest halts, as $h(C_{2})=1$ represents the least halt number as both objects must be worked by both agents;
moreover the exchange must occur at the halfway point of the execution
of the scheme in order for the scheme to be optimal. We may continue
then under the assumption that $n\geq3$.  Due to the role they will play as remainders in the Euclidean algorithm, we shall denote the number of type 1 and type 2 agents by $r_1$ and $r_2$ respectively and without loss we take
$r_1 > r_2$ with $r_1 = ar_2 + r$ for positive integers $a$ and $r$ ($1\leq r\leq r_2 -1$). (We assume that the speeds of the two agent types differ, but do not specify which is the faster.)
An optimal (uniform) scheme is characterised by the
condition that each object is worked by agents of type 1 and type 2 respectively for $r_1 = ar_2 + r$ a.u. and $r_2$ a.u. It follows that $n = (a+1)r_2 + r$.

\vspace{.2cm}

\textbf{Definition 3.1 (Euclidean scheme) }We consider the following
refined algorithm, which has stages corresponding to the lines of
the Euclidean algorithm applied to the pair $(r_1,r_2)$. The associated
optimal scheme will be known as the \emph{Euclidean $(r_{1},r_{2})$-scheme},
denoted by $E=E(r_1,r_2)$.

Let the opening line of the Euclidean algorithm be $r_1 
 = a_1r_2 + r_3$.
We suppress subscripts and write $a$ and $r$ for $a_1$ and $r_3$ respectively. 
Partition the objects into $a+2$ sets as follows. There
are $a+1$ sets $S_{0},S_{1},\cdots,S_{a}$, each of order $r_2$, and a remainder set, $S_{a+1}$, of order $r$. Denote by $S^{(p)}$
the union 
\[
S^{(p)}:=S_{0}\cup S_{1}\cup\cdots\cup S_{p}\,(0\leq p\leq a+1).
\]

At the beginning of Stage 1, each object of $S_0$ is assigned to an agent of type 2, while all other objects are assigned to type 1 agents. 

In Stage 1, all agents first work their initial object for a time of $r_2$ a.u. 

The objects of $S_0$ are then exchanged with the objects of $S_1$ so that each agent that was working on an object of $S_0$ changes to an object from $S_1$ and vice-versa. 

After a further $r_2$ a.u., the objects of $S_1$ exchange agents with those of $S_2$. This action of object exchange between agents is carried out in this fashion so that after the $(i+1)$st halt the objects of $S_i$ are exchanged with those of $S_{i+1}$, until on the $a$th exchange the objects of $S_{a-1}$ and $S_a$ are exchanged between the corresponding agent sets. This entails $a$
halts in all. Stage 1 ends with the $a$th exchange.

At this transition point between Stages 1 and 2, all objects in $S^{(a-1)}$
have been worked by type 2 agents for $r_2$ a.u. and by type 1
agents for $(a-1)r_2$ a.u. In the subsequent stages, which collectively represent a time interval of  $n - ar_2 = (a+1)r_2 + r_3 - ar_2 = r_2 + r_3$ a.u., each of the members
of $S^{(a-1)}$ will continue to be worked by their current type 1 agent. This in effect removes these $ar_2$ type 1 agents from further consideration, leaving an $(r_2 + r_3)$-problem with the type 2 agents now in the majority. Since the objects in $S^{(a-1)}$ undergo no
further agent exchange, we shall say they have become \emph{passive} objects, while the other objects, which now enter the second stage, are labelled \emph{active}. Similarly an agent is described as passive or active according as the agent is working a passive or active object. 

We now repeat the preceding process recursively, mirroring the action of the Euclidean algorithm itself.  Stage 2 acts the process of Stage 1 for the objects of $S_{a}\cup S_{a+1}$
using the second line of the Euclidean algorithm for $(r_1,r_2)$,
but with the roles of the type 1 and type 2 agents reversed. This role alternation of the two types is a feature that  persists as the process is acted throughout subsequent stages.  

Acting the Euclidean algorithm on $(r_1,r_2)$, where $r_1$ and $r_2$ have no common factor, yields for some $t\geq2$: 

\[
r_{1}=a_{1}r_{2}+r_{3}
\]
\[
r_{2}=a_{2}r_{3}+r_{4}
\]
\[
\vdots
\]

\[
r_{i}=a_{i}r_{i+1}+r_{i+2}
\]
\[
\vdots
\]
\[
r_{t-1}=a_{t-1}r_{t}+1
\]
\[
r_{t}=r_{t}\times1+0.
\]

Stage $i$ will correspond to line $i$ of the Euclidean algorithm.  In particular we have $(r_t,a_t,r_{t+1},r_{t+2})=(r_t,r_t,1,0)$.
Stage $t$ then has $a_t = r_t$ exchanges, at the end of which all objects are worked by their current agent for the final a.u., taking all builds to completion and marking the conclusion of what we deem to be Stage $t$. 

\vspace{.2cm}

The sets $S_i$ and the coefficient $a$ change as we pass from one stage of the algorithm to the next. To compare these sets between stages calls for notation with double subscripts. 

To compare the sets of objects in question, consider Stage $i$ $(1\leq i\leq t)$ of $E(r_1,r_2)$. Let $S_{i,j}$ denote the set $S_j$ at the $i$th stage of the scheme $(0\leq j\leq a_{i}+1)$. By construction, the active objects of Stage $i$ all begin with identical work records with respect to both agent types. Moreover, during the $i$th stage, each object in $S_{i,j}$ has been worked by each type of agent for the same length of time as each other member of $S_{i,j}$; indeed we say something more precise.

\vspace{.2cm}

\textbf{Definition 3.2 }The \emph{work record} for the objects in $S_{i,j}$, is $(t_{1},t_{2})_{i,j}$, where $t_{p}$ $(p=1,2)$ denotes the number of a.u. for which the members of $S_{i,j}$ have been worked by agents of type $p$ at the completion of Stage $i$.

\vspace{.2cm}

\textbf{Remark 3.3} It follows from the recursive construction of the stages of the algorithm that for any $i$ there are only two distinct work record pairs determined by whether, at the  end of Stage $i$, $S_{i,j}$  represents a set of passive objects $(0\leq j\leq a_i-1)$ or active objects $(a_i \leq j\leq a_i+1)$. We therefore simplify notation by denoting the two respective pairs by $(t_1,t_2)_i$ and $(T_1,T_2)_i$. In particular all objects active in Stage $i + 1$ begin that stage having been worked $T_1$ a.u. by type 1 agents and $T_2$ a.u. by type 2 agents.

\vspace{.2cm}

The following proof of optimality of $E(r_1,r_2)$ includes the precise work record of all objects that become passive at the conclusion of each stage.

\vspace{.2cm}

\textbf{Theorem 3.4} The Euclidean scheme $E(r_1,r_2)$ is optimal. 

\begin{proof}

We show that in $E(r_1,r_2)$ each object is worked continuously and uniformly for a total of $H = r_1 + r_2 = n$ a.u, from which the result follows from Proposition 2.1.3(b). The proof is by induction on the stage number $i$.  

An object that becomes passive at the end of Stage 1 has work record: 
\begin{equation}
(t_1,t_2)_1 = ((a_1-1)r_2,r_2) = (a_1r_2 - r_2, r_2) 
= (r_1 - r_2 - r_3, r_2).
\end{equation}
Similarly at the close of Stage 1 we have  for the active objects 
\begin{equation}
 (T_1,T_2)_1 = (a_1r_2,0) = (r_1 - r_3, r_2 - r_2). 
\end{equation}

The length of Stage 1 is $a_1r_2 = r_1 - r_3$. The remaining time before $H$ expires matches the number of active objects and agents as all three equal

\[n - (r_1 - r_3) = r_1 + r_2 - r_1 + r_3 = r_2 + r_3.\]

Therefore in moving from Stage 1 to Stage 2, the framework for the scheme is repeated with $(r_1,r_2)$ replaced by $(r_2,r_3)$, where $r_2 = a_2r_3 + r_4$ (so that $r_2 > r_3)$. However the roles of type 1 and type 2 agents are reversed, as it is type 2 that now form the majority. 

An object that becomes passive at the end of Stage 2 has work record:
\begin{equation}
(t_1,t_2)_2 = (r_1 - r_3,0) + (r_3, (a_2 - 1)r_3) =
(r_1 - r_3,0) + (r_3, r_2 -r_3 -r_4)
= (r_1, r_2 - r_3 - r_4).
\end{equation}

Similarly, at the close of Stage 2 we have  for the active objects
\begin{equation}
(T_1,T_2)_2 = (r_1 - r_3,0) + (0,a_2r_3) = 
(r_1 - r_3, r_2 - r_4). 
\end{equation}

The length of Stage 2 is $a_2r_3 = r_2 - r_4$ and hence the remaining time before $H$ expires is
\[r_2 + r_3 - (r_2 - r_4) = r_3 + r_4, \] which matches the number of active objects and active agents as we enter Stage 3. 

We next show by induction on $i$ that for the objects that become passive at the end of Stage $i$, if $i$ is odd, then 
\begin{equation}
(t_1,t_2)_i = (r_1 - r_{i+1} - r_{i+2},r_2) 
\end{equation}
while for the active objects we have
\begin{equation}
(T_1,T_2)_i = (r_1 - r_{i+2},r_2 - r_{i+1}).
\end{equation}
If $i$ is even then
\begin{equation}
(t_1,t_2)_i = (r_1,r_2 - r_{i+1} - r_{i+2}) 
\end{equation}
\begin{equation}
(T_1,T_2)_i = (r_1 - r_{i+1},r_2 - r_{i+2}).
\end{equation}
Moreover, at the conclusion of Stage $i$, the number of active objects matches both the number of active agents and the time remaining until the harmonic optimum expires, which is $r_i + r_{i+1}$. The induction is anchored on the $i = 1$ and $i = 2$ cases given in (3) - (6). 

Let $i+1 \geq 3$ be odd.  Stage $i+1$ is based on 
\[r_{i+1} = a_{i+1}r_{i+2} + r_{i+3}. \] Since $i$ is even we apply (10) to give $(T_1,T_2)_i = (r_1 - r_{i+1}, r_2 - r_{i+2})$ for the first term in the next sum.  For the second term, we apply the form of  $(t_1,t_2)_1$ given in (3) but increment the subscripts by $i$. This yields: 

\[(t_1,t_2)_{i+1} = (r_1 - r_{i+1}, r_2 - r_{i+2}) + (r_{i+1} - r_{i+2} -  r_{i+3}, r_{i+2})   \]
\[= (r_1 - r_{i+2} - r_{i+3}, r_2) \] in accord with (7) for $i+1$; similarly adjusting (4) for line $i+1$ yields 
\[(T_1,T_2)_{i+1} = (r_1 - r_{i+1},r_2 - r_{i+2}) + 
(r_{i+1} - r_{i+3},0) = (r_1 - r_{i+3}, r_2 - r_{i+2})\] 
 in accord with (8) for $i+1$.

Now let $i+1\geq 4$ be even. Since $i$ is odd we apply (8) to give the first term of the next sum, which is
$(T_1,T_2)_i = (r_1 - r_{i+2}, r_2 - r_{i+1})$,
while for the second term we apply the form of  $(t_1,t_2)_2$ in (5) but increment the subscripts of $(r_3,r_2 - r_3 - r_4)$ by $i - 1$ (as $2 + (i - 1) = i + 1$). This yields: 
\[(t_1,t_2)_{i+1} = (r_1 - r_{i+2}, r_2 - r_{i+1}) +
(r_{i+2}, r_{i+1} - r_{i+2} - r_{i+3} ) = \]
\[(r_1 , r_2 - r_{i+2} - r_{i+3})\] in accord with (9) for $i+1$; similarly using (6) we infer
\[(T_1,T_2)_{i+1} = (r_1 - r_{i+2}, r_2 - r_{i+1}) + 
(0, r_{i+1} - r_{i+3}) = \]
\[(r_1 - r_{i+2}, r_2 - r_{i+3}) \]
in accord with (10) for $i+1$, and so the induction continues in both cases. 

At the transition from Stage $i+1$ to Stage $i+2$, for both the odd and even case, the common number of active objects and agents is  
\[r_{i+1} + r_{i+2} - a_{i+1}r_{i+2} = \]
\[r_{i+1} +r_{i+2} - (r_{i+1} - r_{i+3}) = 
r_{i+2} + r_{i+3},\] which equals the remaining time in the harmonic optimum:
$r_{i+1} + r_{i+2} - a_{i+2}r_{i+2}$, and so the induction continues. 

Putting $i = t$ in (7) gives that at the conclusion of the final stage, if $t$ is odd, 
\[(t_1,t_2)_t = (r_1 - r_{t+1} - r_{t+2},r_2) = \]
\[(r_1 - 1 - 0, r_2) = (r_1 - 1,r_2). \] 
The final a.u. then adds $(1,0)$ to this pair for all passive objects, giving the required final pair of $(r_1,r_2)$ representing the numbers of a.u. worked by type 1 and type 2 agents respectively. On the other hand by (8) we get
\[(T_1,T_2)_t = (r_1 - r_{t+2},r_2 - r_{t+1}) = (r_1,r_2 -1),\] and this time it is $(0,1)$ added to work record of the final active object, giving the required pair, $(r_1,r_2)$.  

If $t$ is even then by (9) we get
\[ (t_1,t_2)_t  =
(r_1, r_2 - r_{t+1} - r_{t+2}) = \]
\[ (r_1, r_2 - 1 - 0) = (r_1,r_2 - 1) .\] The final a.u. adds $(0,1)$ to this pair for all passive objects, giving the required final pair of $(r_1,r_2)$.   On the other hand by (10) we get
\[(T_1,T_2)_t = (r_1 - r_{t+1},r_2 - r_{t+2}) = (r_1 - 1,r_2).\] The single active object this applies to then has $(1,0)$ added to complete the scheme with the required pair of $(r_1,r_2)$.  

In all cases then, each object is worked for $r_1 + r_2 = n$ a.u $= H$, and so $E(r_1,r_2)$ is optimal.  

\end{proof}

We record some useful observations that emerged in the previous proof.

\vspace{.2cm}

\textbf{Corollary 3.5} For the Euclidean scheme $E=E(r_{1},r_{2})$:

(a) there are $t$ stages in $E$. The halt number $h=h(E)$ is $h=\sum_{i=1}^{t}a_{i}$; 

(b) The length of Stage $i$ is $a_ir_{i+1}=r_{i} -r_{i+2}$ a.u.
$(1\leq i\leq t-1)$, and Stage $t$ has length $a_t+1$ a.u.; 

(c) In any stage, the members of $S_0\cup S_a$ undergo a single exchange, those of $\cup_{i=1}^{a-1}S_i$ undergo two exchanges, while those of $S_{a+1}$ are not exchanged. 

\begin{proof}
 
(a) For Stage $i$ of the $t$ stages there are $a_i$ halts including for Stage $t$. Summing the $a_i$ thus gives $h(E)$. 

(b) For $1\leq i\leq t$, for Stage $i$, prior  to each of the $a_i$ halts, all agents work continuously for $r_{i+1}$ a.u. Hence the length of Stage $i$ is $a_{i}r_{i+1}$ a.u. $(1\leq i\leq t-1)$ and Stage $t$ has $a_tr_{t+1} + 1 = a_t + 1$ a.u.  

(c) These are observations of the exchange process.

\end{proof}

\vspace{.2cm}

\textbf{Example 3.6} $r_1 = 180, r_2 = 53,n = r_1 +r_2 = 233.$

\[
180=3\times53+21
\]
\[
53=2\times21+11
\]
\[
21=1\times11+10
\]
\[
11=1\times10+1
\]
\[
10=10\times1+0.
\]

We shall initially label our parameters $r_1, r_2, a,r$ where $r_1 = ar_2 + r$, updating values for these parameters in accord with our scheme as we pass from one stage to the next. 

Set the parameters for Stage 1 based on the first line of the Euclidean
algorithm for $(r_1, r_2)$: 

\[r_1 = 180,r_2 = 53, a=3,r=21, n = r_1 + r_2  = 233.\]
Since $a=3$, we partition the objects into $a+2=5$ sets: 
\[
S^{(a-1)}=S^{(2)}=S_0\cup S_{1}\cup S_2,\,|S_{i}|= r_2 =53\,(0\leq i\leq3),
\]
\[|S_4|=r=21. \]

All agents now work their objects for $r_2 = 53$ a.u. and then halt to allow the first exchange. There are $a=3$ such halts and exchanges in Stage 1. 

At the end of Stage $1$ we have: 
\[
(t_1,t_2)_1=((a-1)r_2,r_2) = (106,53),
\]
\[(T_1,T_2)_1 = (ar_2,0)=(159,0).\]
Each object in $S^{(a-1)}=S^{(2)}$ will next be worked on to completion by the type 1 agent that has just taken it into possession, making $ar_2 = 3r_2$ type 1 agents unavailable for the objects of $S_3\cup S_4$.
The number of a.u. until completion is given by $r_2 + r.$
Therefore at the end of the entire process the members of the sets in $S^{(2)}$ will have the required final assignment pair: 
\[
((a-1)r_2 + r_2 + r,r_2)=(ar_2 + r,r_2)=(r_1,r_2),
\]
which in this case this gives:
\[
(106 + 53 + 21,53) = (180,53) = (r_1,r_2).
\]
There remain $r$ type 1 and $r_2$ type 2 agents available to $S_3\cup S_4$ at the beginning of the second stage. 

We now analyse Stage 2 for the objects in $S_3\cup S_4$. As we pass from one stage to the next, the roles of the agent types are reversed. Our parameters have their values updated. Note that the entries of $(T_1,T_2)_1$ are the respective inherited starting values of type 1 a.u. and type 2 a.u. for the remaining $r_2 + r$ active objects of Stage 2. 

\vspace{.2cm}

Stage 2: 
\[r_1 = 53,r_2 = 21, a=2, r=11,\,n = r_1 + r_2 = 74.\]
\[ S^{(a-1)}=S^{(1)}=S_0\cup S_1,\,|S_i|= r_2 =21,\,(0\leq i\leq2), \]
\[|S_3|=r = 11.\]
\[ (t_1,t_2)_2=(159,0)+(r_2,(a-1)r_2)=(159,0)+(21,21)=(180,21), \]
\[ (T_1,T_2)_2=(159,0)+(0,ar_2)=(159,0)+(0,42)=(159,42).\]

Note that the form of the pair $((a-1)r_2,r_2)$ of Stage 1 is reversed in Stage 2 to become $(r_2,(a-1)r_2)$; this is due to the exchange of roles of type 1 and type 2 agents, a feature of each passage to a new stage. 

Since $r_2 + r= 21 + 11 = 32,$ at the end of the entire process, the objects in the current $S^{(a-1)}=S^{(1)}$ will have final assignment pair $(180,21+32)=(180,53)$, as required. 
\vspace{.2cm}

Stage 3:
\[r_1 = 21,\,r_2 = 11,\,a=1,\,r=10,\,n = r_1 + r_2 = 32.\]
\[ S^{(a-1)}=S_0,\,|S_i|=r_2 = 11,(0\leq i\leq1),|S_2|=r=10. \]
\[
(t_1,t_2)_3 = (159,42)+((a-1)r_2,r_2) = (159,42) + (0,11) = (159,53),
\]
\[ (T_1,T_2)_3 = (159,42)+(ar_2,0)=(159,42)+(11,0)=(170,42).\]

Since $r_2 + r = 11 + 10 = 21$, at the end of the entire process, the objects in the current $S^{(a-1)}=S^{0}$ will have final assignment pair $(159+21,53)=(180,53)$, as required. 

\vspace{.2cm}

Stage 4: 
\[r_1 = 11,\,r_2 = 10,\,a=1,\,r=1,\,n= r_1 + r_2 = 21.\]
\[S^{(a-1)}=S_0,\,|S_i|=r_2 = 10,(0\leq i\leq1),\,|S_2| = r = 1.\]
\[
(t_1,t_2)_4=(170,42)+(r_2,(a-1)r_2)=(170,42)+(10,0)=(180,42),
\]
\[ (T_1,T_2)_4=(170,42)+(0,ar_2)=(170,42)+(0,10)=(170,52).\]

Since $r_2 + r = 10 + 1 = 11$, at the end of the entire process, the objects in the current $S^{(0)}$ will have final assignment pair $(180,42+11)=(180,53)$, as required. 

\vspace{.2cm}

Stage 5 (final stage)
\[r_1 = 10,\,r_2 = 1,\,a = 10,\,r=0,\,n = r_1 + r_2 = 11.\]
\[S^{(a-1)}=S_0\cup S_{1}\cup\cdots\cup S_9,\,|S_i|= r_2 = 1,\,(0\leq i\leq10),\,|S_{11}| = r = 0.\]
\[
(t_1,t_2)_5 = (170,52)+((a-1)r_2 + 1,r_2)=(170,52)+(9+1,1)=(180,53);
\]
\[ (T_1,T_2)_5 = (170,52) + (ar_2,1) = (170,52)+(10,1)=(180,53).\]

By Proposition 3.2(a), the total number of halts $h$ of the agents is the sum of values of the parameter $a$ over the five stages:
\[
h(180,53) = 3 + 2 + 1 + 1 + 10 = 17.
\]
This compares with the simple cyclic method $C_{233}$, which has
$n - 1 = 232$ such pauses in production. 

Proposition 3.2(b) gives the list of lengths of the five stages as,
$159,42,11,10$ and $11$, which sum to $n = 233$. 

\section{Euclidean schemes and their halt numbers}

\subsection{Matrix of an optimal scheme}

With each optimal $n$-scheme $S$ of the $m=2$ problem we may associate an
$n\times n$ binary matrix $M=M(S)$ whereby the $(i,j)$th entry
of $M$ is $0$ or $1$ according as object $O_i$ in the $j$th
time interval $I_j$ is worked by an agent of type 1 or type 2. 

The labelling of the set of objects is arbitrary: a permutation $\pi$
of the rows of $M$ gives a new matrix $M_{\pi}$ that corresponds
to a permutation of the object set ${\cal O}$, which is to say a
re-numbering of the members of ${\cal O}$ (by $\pi^{-1}$). Hence the corresponding schemes, $S(M)$ and $S(M_{\pi})$, are equivalent up to the labelling of the
objects. 

At the beginning of each interval in the execution of an optimal scheme $S$,
the type 1 and type 2 agents are re-assigned their objects, forming
two sets $K_1$ and $K_2$ of objects respectively. For a given
ordering of ${\cal O}$, these sets are equal for all stages for a
pair of schemes $S$ and $S'$ exactly when $M(S)=M(S')$. Therefore
we may declare that two schemes $S$ and $S'$ are \emph{equivalent}
if $M(S)=M(S')$, and then extend this notion of equivalence to the
case where $M(S)=M_{\pi}(S')$, for some row permutation $\pi$ of
the rows of $M(S')$. 

Clearly the columns of $M(S)$ will contain $k = r_1$ entries of
$0$ and $r_2$ entries of $1$. Proposition 2.1.5(c) is equivalent
to the statement that $S$ is optimal if and only if each row also
contains exactly $k$ instances of $0$. 

In summary, $S$ is optimal if and only if $M(S)$ is $k$\emph{-uniform},
meaning that each row and column of $M$ has exactly $k$ zeros. Up
to equivalence, there is a one-to-one correspondence between optimal
$n$-schemes $S$ with $k$ type 1 agents and $k$-uniform $n\times n$
binary matrices. 

\subsection{Successive Fibonacci numbers}

The number of halts in a Euclidean scheme is determined by the sum 
of the coefficients $a_i$, which corresponds to the interpretation
of the Euclidean algorithm whereby each step involves simply subtracting
the smaller of the two integers in hand from the larger. In the case
of a pair of successive Fibonacci integers, $f_{p+1},f_p$ say,
all the $a_i$ are equal to $1$ (apart from the final coefficient
where the remainder is $0$). This leads to a halt number of order
$\log n$, where $n=f_{p+1}+f_p = f_{p+2}$. 

To see this put $r_1=f_{p+1}$, $r_2 = f_p$, two consecutive
and distinct Fibonacci numbers (so that $p\geq2$), whence $n=f_{p+1}+f_p=f_{p+2}$.
There are then $p-1$ lines in the Euclidean algorithm:
\[
f_{p+1}=f_p+f_{p-1}
\]
\[
f_p = f_{p-1}+f_{p-2}
\]
\[
\vdots
\]

\[
f_4 = f_3 + f_2 \Leftrightarrow3=2+1
\]
\[
f_3=2f_2 + 0\Leftrightarrow2=2\cdot1+0.
\]

Now $a_i = 1$ for all but the final line where the coefficient of
$f_2$ equals $2$. By Corollary 3.5(a), $h(E(f_{p+1},f_p))=(p-1-1)+2=p.$
The respective lengths of the $p-1$ stages are $f_p,f_{p-1},\cdots,f_3 = 2,2f_2 + 1 = 3$.
Since $f_2 = f_1$, the sum of these lengths is equal to 
\[
(f_p + f_{p-1}+\cdots + f_2 + f_1 ) + 1=(f_{p+2}-1)+1=f_{p+2} = n, \]
all in accord with Corollary 3.5(b). Since $f_n=\lfloor\frac{\phi^{n}}{\sqrt{5}}\rceil$,
the nearest integer to $\frac{\phi^{n}}{\sqrt{5}}$, where $\phi$
denotes the Golden ratio $\frac{1+\sqrt{5}}{2}$, (see, for example,
{[}5{]}), it follows that for all sufficiently large $p$, 
\[
n=f_{p+2}=\lfloor\frac{\phi^{p+2}}{\sqrt{5}}\rceil\geq\phi^{p}.
\]
Hence, since $h=p$ we then have 
\[
n\geq\phi^{h}\Rightarrow h\leq\log_{\phi}n=\frac{\ln n}{\ln\phi}\approx2.078\ln n.
\]
\[
\Rightarrow  h(E(f_p,f_{p+1}))=O(\log f_{p+2})=O(\log n).
\]

\newpage

\textbf{Example 4.2.1 }We represent the Euclidean Fibonacci scheme
$F_p$, for $p=5$ in matrix form. Since $f_5=5$ and $f_6=8$,
there are $5$ agents labelled type 2 ($A_1$ -$A_5$) and $8$ of type 1 ($A_6$ -$A_{13}$), with $p = 5$ stages, and $f_7 = 13$ objects in all. By Corollary 3.5(b), the halts, which number $p=5$, occur
at the end of the intervals $f_5 = 5$, $5+f_4 = 8$, $8+f_3 = 10$,
$10 + f_2 = 11$, and $11+f_1 = 12$. 

The matrix $M(F_5)$ of the Euclidean Fibonacci $5-$scheme is shown
in Figure 2. In all but the final stage the active objects form three sets, $S_0,S_1,$ and $S_2$, which feature only one set of exchanges, which occur between $S_0$ and $S_1$. Initially
object $O_i$ is worked by agent $A_i$.  The columns within each stage are identical. As we pass from one stage to the next, the objects of $S_0$ become
passive, the `old' $S_2$ becomes the `new' $S_0$, while the
`old' $S_1$ splits to form the `new' $S_1,S_2$ pair. The subscript on the $(i,j)$th entry gives
the number of the agent that is working object $O_i$ in the $j$th a.u. of the scheme $F_5$ (and so, in respect to $M(F_5)$, the objects label the rows while the agents, recorded as subscripts, may `move' from one column to the next).

\vspace{.2cm}

Stage 1: 
\[
r_1 = 8,r_2 = 5,a=1,r=3,n = r_1 + r_2 = 13.
\]
\[
S^{(a-1)}=S_0,|S_i|= r_2 = 5,(0\leq i\leq1),|S_2|=r=3.
\]
\[
S_0=\{O_1,O_2,\cdots,O_5\},S_1=\{O_6,O_7,\cdots,O_{10}\},S_2=\{O_{11,}O_{12},O_{13}\}.
\]
At the conclusion of Stage 1, the objects in $S_0$ are with their
set of final agents, $\{A_6,A_7,A_8,A_9,A_{10}\}$, (see Fig.
2.)

\vspace{.2cm}

Stage 2: 
\[
r_1 = 5, r_2 = 3,a=1,r=2,n= r_1 + r_2 =8.
\]
\[
S^{(a-1)}=S_0,|S_i|= r_2 =3,(0\leq i\leq1),|S_2|=r=2.
\]
\[
S_0= \{O_{11},O_{12},O_{13}\},S_1= \{O_8,O_9,O_{10}\},S_2=\{O_6,O_7\}.
\]

At the end of Stage 2, the objects in $S_{0}$ are with their
final agents, \newline $\{A_3,A_4,A_5\}$, (see Figure 2). 

\vspace{.2cm}

Stage 3: 
\[
r_1 = 3, r_2 = 2,a=1,r=1,n= r_1 + r_2 = 5.
\]
\[
S^{(a-1)}=S_0,|S_i|= r_2 =2,(0\leq i\leq1),|S_2|=r=1.
\]
\[
S_0=\{O_6,O_7\},S_1=\{O_8,O_9\},S_2=\{O_{10}\}.
\]

At the end of Stage 3, the objects in $S_0$ are with their
final agents, $\{A_{11},A_{12}\}$. 

\vspace{.2cm}

Stage 4: 
\[
r_1 = 2,r_2 = 1 ,a=2,r=0,n = r_1 + r_2 = 3.
\]
\[
S^{(a-1)}=S_{0}\cup S_{1},|S_{i}|= r_2 =1,(0\leq i\leq2).
\]
\[
S_{0}=\{O_{10}\},S_{1}=\{O_{9}\},S_{2}=\{O_{8}\}.
\]

At the end of Stage 4, the object in $S_0$ is with its final
agent, $A_2$.

\vspace{.2cm}

Stage 5: 
\[
r_2 = r_1 = 1,a=1,r=0,n = r_1 + r_2 = 2.
\]
\[
S^{(a-1)}=S_{0},|S_{i}|= r_2 =1,(0\leq i\leq1).
\]
\[
S_0=\{O_8\},S_1=\{O_9\}.
\]

At the conclusion of Stage 5, the objects in $S_0$ and $S_1$
are with their respective final agents, $A_{13}$ and $A_1$. 

\vspace{.3cm}

\begin{tabular}{|c|c|c|c|c|c|c|c|c|c|c|c|c|c|}
\hline 
 &  &  &  &  & 5 &  &  & 8 &  & 10 & 11 & 12 & \tabularnewline
\hline 
\hline 
$O_{1}$ & $0_{1}$ & $0_{1}$ & $0_{1}$ & $0_{1}$ & $0_{1}$ & $1_{6}$ & $1_{6}$ & $1_{6}$ & $1_{6}$ & $1_{6}$ & $1_{6}$ & $1_{6}$ & $1_{6}$\tabularnewline
\hline 
$O_{2}$ & $0_{2}$ & $0_{2}$ & $0_{2}$ & $0_{2}$ & $0_{2}$ & $1_{7}$ & $1_{7}$ & $1_{7}$ & $1_{7}$ & $1_{7}$ & $1_{7}$ & $1_{7}$ & $1_{7}$\tabularnewline
\hline 
$O_{3}$ & $0_{3}$ & $0_{3}$ & $0_{3}$ & $0_{3}$ & $0_{3}$ & $1_{8}$ & $1_{8}$ & $1_{8}$ & $1_{8}$ & $1_{8}$ & $1_{8}$ & $1_{8}$ & $1_{8}$\tabularnewline
\hline 
$O_{4}$ & $0_{4}$ & $0_{4}$ & $0_{4}$ & $0_{4}$ & $0_{4}$ & $1_{9}$ & $1_{9}$ & $1_{9}$ & $1_{9}$ & $1_{9}$ & $1_{9}$ & $1_{9}$ & $1_{9}$\tabularnewline
\hline 
$O_{5}$ & $0_{5}$ & $0_{5}$ & $0_{5}$ & $0_{5}$ & $0_{5}$ & $1_{10}$ & $1_{10}$ & $1_{10}$ & $1_{10}$ & $1_{10}$ & $1_{10}$ & $1_{10}$ & $1_{10}$\tabularnewline
\hline 
$O_{6}$ & $1_{6}$ & $1_{6}$ & $1_{6}$ & $1_{6}$ & $1_{6}$ & $0_{1}$ & $0_{1}$ & $0_{1}$ & $0_{1}$ & $0_{1}$ & $1_{11}$ & $1_{11}$ & $1_{11}$\tabularnewline
\hline 
$O_{7}$ & $1_{7}$ & $1_{7}$ & $1_{7}$ & $1_{7}$ & $1_{7}$ & $0_{2}$ & $0_{2}$ & $0_{2}$ & $0_{2}$ & $0_{2}$ & $1_{12}$ & $1_{12}$ & $1_{12}$\tabularnewline
\hline 
$O_{8}$ & $1_{8}$ & $1_{8}$ & $1_{8}$ & $1_{8}$ & $1_{8}$ & $0_{3}$ & $0_{3}$ & $0_{3}$ & $1_{11}$ & $1_{11}$ & $0_{1}$ & $0_{1}$ & $1_{13}$\tabularnewline
\hline 
$O_{9}$ & $1_{9}$ & $1_{9}$ & $1_{9}$ & $1_{9}$ & $1_{9}$ & $0_{4}$ & $0_{4}$ & $0_{4}$ & $1_{12}$ & $1_{12}$ & $0_{2}$ & $1_{13}$ & $0_{1}$\tabularnewline
\hline 
$O_{10}$ & $1_{10}$ & $1_{10}$ & $1_{10}$ & $1_{10}$ & $1_{10}$ & $0_{5}$ & $0_{5}$ & $0_{5}$ & $1_{13}$ & $1_{13}$ & $1_{13}$ & $0_{2}$ & $0_{2}$\tabularnewline
\hline 
$O_{11}$ & $1_{11}$ & $1_{11}$ & $1_{11}$ & $1_{11}$ & $1_{11}$ & $1_{11}$ & $1_{11}$ & $1_{11}$ & $0_{3}$ & $0_{3}$ & $0_{3}$ & $0_{3}$ & $0_{3}$\tabularnewline
\hline 
$O_{12}$ & $1_{12}$ & $1_{12}$ & $1_{12}$ & $1_{12}$ & $1_{12}$ & $1_{12}$ & $1_{12}$ & $1_{12}$ & $0_{4}$ & $0_{4}$ & $0_{4}$ & $0_{4}$ & $0_{4}$\tabularnewline
\hline 
$O_{13}$ & $1_{13}$ & $1_{13}$ & $1_{13}$ & $1_{13}$ & $1_{13}$ & $1_{13}$ & $1_{13}$ & $1_{13}$ & $0_{5}$ & $0_{5}$ & $0_{5}$ & $0_{5}$ & $0_{5}$\tabularnewline
\hline 
\end{tabular}

\vspace{.3cm}

\centerline{Figure 2: matrix for the Fibonacci scheme $F_5$.}

\vspace{.2cm}

\textbf{Remark 4.2.2 }\emph{Halting number of generic Euclidean examples}.

\vspace{.3cm}

The halting number of the Euclidean scheme based on two successive
Fibonacci numbers that sum to $n$ is of order $O(\ln n)$, while
at the other extreme the only optimal scheme $S$ when $r_1 = n$ and $r_2 = 1$ has a halting number of $h(S)  = n-1$. 

For a given positive integer $n$, the expected halting number is
the mean number of subtractions in the Euclidean algorithm for two
relatively prime integers $r_2 < r_1$ such that $r_1 + r_2 = n.$
The mean number of lines in the Euclidean algorithm in this case is
certainly of order $O(\ln n)$, and indeed much more precise statements
are known (see {[}4{]}, which cites {[}3{]}). General considerations
suggest that the mean halt number for such a pair $(r_1,r_2)$
is then $O((\ln n)^{2})$ [2], but this has not been proved. 

The question may be formulated as asking for the expected sum of the
coefficients of the continued fraction expansion of $\frac{a}{n-a}$
for $a$ drawn at random from $1,\cdots,\frac{n}{2}$ (and restricted
to be coprime to $n$). From the Gauss-Kuzmin distribution {[}6{]}
and the St Petersburg paradox analysis, Terry Tao comments in {[}2{]}
that it is reasonable to expect the answer to be $O((\ln n)^{2})$,
though rigorous analysis may be difficult. 

\subsection{Production times with allowance for handovers}

Consider again the $m = 2$ case with timings $t_1 = 1$ for a set of $r_1$ agents, and $t_2 = T \neq 1$, for a set of $r_2$ agents. The optimal production time is:
\begin{equation}
H=n(r_1 + \frac{r_2}{T})^{-1}=\frac{nT}{r_1T +r_2}.
\end{equation}
We now make allowance for a time interval, the \emph{halting time},
of length $\varepsilon>0$ for each halt during a scheme and an initial
\emph{loading time} of $\varepsilon$ also. For $C_{n}$ the total
production time $C$ is then:
\begin{equation}
C=\frac{nT}{r_1T + r_2 } + n\varepsilon.
\end{equation}

Applying (12) to Example 4.2.1, taking $T = 2$ and $\varepsilon=0.005$
we have to five significant figures that the harmonic optimum time
is:
\[
H=\frac{13\times2}{8\times 2 + 5} = \frac{26}{21} = 1.2381.
\]
This compares to the production times $C$ of $C_{13}$ and $E$ of the Euclidean scheme: 
\[
C=H+13\varepsilon=1.3031,\,\,E=H+6\varepsilon=1.2681.
\]
The excess percentages over the harmonic optimum are respectively $5.3\%$ and $2.4\%$.

Applying (12) to Example 3.6 however, again taking $T=2$ and $\varepsilon=0.005$, to five significant figures we get:
\[
H=\frac{233\times2}{180\times2 + 53} = 1.1283;\,\,C = H + 233\varepsilon= 2.2933.
\]
Hence with an exchange period $\varepsilon$ equal to 18 seconds we
have for $C_{233}$ an increase in production time above the harmonic optimum of $103$\%. In contrast, with the Euclidean scheme we calculated $h(E(53,180))=17$, so that the build time $E$ of $E(53,180)$
is:
\[
E=H+18\varepsilon=H + 0.09 = 1.2183,
\]
which represents an increase of only 8\% above the harmonic optimum. This suggests that with a generic example involving non-zero exchange
times, the Euclidean scheme is much more efficient that the $n$-cyclic scheme. 

\subsection{Comparison with Biker-hiker problem}

As mentioned in the introduction, in the $m = 2$ case the $n$-object problem with $k$ machines corresponds to the Biker-hiker problem with $k$ bicycles and $n$ travellers where the bicycles may move instantaneously from their current staging post to any other.  It follows that the time of an optimal scheme where this superpower may be exploited must be less than the least time possible in the original setting. 

It transpires that for the Biker-hiker problem, in any optimal solution each traveller covers exactly $\frac{k}{n}$ of the \emph{distance} of the journey on the faster mode of transport (the bicycle), while in our $n$-object problem with $k$ faster agents, each object spends exactly $\frac{k}{n}$ of the total journey \emph{time} with faster agents (remembering that when $m = 2$, all optimal schemes are uniform).  In the Biker-hiker case, the time taken cycling will be less than $\frac{k}{n}$ of the total scheme time just because cycling beats walking, and so will cover the $\frac{k}{n}$ distance in less than that proportion of the total time.

This is confirmed by comparing the formulas for the optimal times.  From (11) we see that if there are $r_1$ agents of one completion time, which we take to be 1, and $r_2$ of another time $T$, then the optimal time is

$$\frac{nT}{r_1T + r_2}.$$

This holds whether or not the first agent type is the faster of the two.  The time required to execute an optimal scheme for the Biker-hiker problem is then:

$$\frac{r_1 + r_2T}{n}.$$

(It may be noted that the lesser time of $1$ and $T$ must apply to bicycle travel. The opposite scenario corresponds to broken bicycles that impede the progress of a traveller obliged to walk them along.  In this case, the optimal completion time is simply the time taken to complete the journey with a broken bicycle, and this is the same for any positive value of $k$, the number of bikes.)  

We then expect the first of these two expressions to be less than the second.  Simplifying this inequality shows it to be equivalent to $(T-1)^2 > 0$, which is true as $T \neq 1$.

\subsection{Euclidean scheme is greedy}

The Euclidean scheme represents a Greedy algorithm in that the process
is never halted unless continuation would result in some object exceeding
its quota of agent type for an (optimal) uniform scheme. 

We show that if the Greedy principle is adopted to give a scheme $S$, we are effectively forced into the Euclidean scheme. After $r_2$ a.u. the set $S_0$ of objects being worked by the members of the
set of type $2$ agents reach their type 2 quota, and so the process
must halt. The members of $S_0$ are then exchanged with a subset $T_0$ of the type 1 agents. The Euclidean directive however is not the only possible continuation of our Greedy scheme $S$. 

Let ${\cal O}$ and ${\cal A}$ denote the $n$-sets of objects and
agents. In general, a Greedy scheme $S$ acts a bijection $\pi:{\cal O}\rightarrow{\cal A}$
subject to the constraint that $S_{0}\pi\cap A_0=\emptyset$, where
$A_0$ denotes the initial set of agents working the members of
$S_0$. For $E(r_1,r_2)$ this bijection $\pi$ has the additional property that $S_0\pi^{2} = S_0$. 

However, in both $E(r_1,r_2)$ and $S$, there is a set of $r_2$
objects that have completed their type 2 quota and which are now assigned
to type 1 agents, there is a set of $r_2$ objects that have been
worked for $r_2$ a.u. by type 1 agents and which now are assigned
to type 2 agents, while the remainder of the objects have been worked
by type 1 agents for $r_2$ a.u. and remain assigned to type 1 agents.
The outcome from both schemes is therefore identical up to the numbering
of the agents in that there is a permutation $'$ of ${\cal O}$ such
that for each $O\in{\cal O}$, at any time point in the execution
of the schemes $E$ and $S$, $O$ in scheme $E$ and $O'$ in scheme
$S$ have been worked for the same length of time by type 2 agents
(and hence also by type 1 agents). Moreover, this equivalence will
persist as we pass to subsequent exchanges and stages provided we
adhere to the Greedy principle of never halting until forced to in
order not to exceed type quotas of objects. Therefore we may identify
any scheme based on the Greedy principle with the Euclidean scheme
$E(r_1,r_2).$ 

It remains to be determined however if an approach based on the Greedy
principle is effective for the general $n$-object problem when $m\geq3$. 

\vspace{.4cm}

ACKNOWLEDGEMENT The thoughtful and astute suggestions of the referees were helpful and appreciated.

\end{document}